%% file: maxwell_time-periodic.tex
\documentclass[12pt,a4paper,twoside]{article}

\title{\sc Theoretical Considerations on the Computation of Generalized Time-Periodic Waves}
\def\shorttitle{Computation of Generalized Time-Periodic Waves}
\def\pauthor{Dirk Pauly and Tuomo Rossi}

\def\mylabelonoff{off}
\def\allowdisbrk{no}

\input head_paper

\input pauly_kopf

\renewcommand{\rot}{\pd\,}
\renewcommand{\Rot}{\rot\!{}_{\dom}}
\newcommand{\Rotb}{\rot\!{}_{\domb}}
\renewcommand{\ROT}{\rot\!{}_{\circ}}
\renewcommand{\pdiv}{\delta}
\renewcommand{\Div}{\pdiv_{\dom}}

\renewcommand{\DIV}{\pdiv}
\renewcommand{\circledast}{*_{\dom}}
\renewcommand{\pr}[4]{{}_{#3}\overset{#4}{\mathsf{D}}{}^{#1}_{#2}}
\renewcommand{\pdi}[4]{{}_{#3}\overset{#4}{\mathsf{\Delta}}{}^{#1}_{#2}}
\renewcommand{\lz}{\mathsf{L}^2}
\renewcommand{\qlz}[1]{\mathsf{L}^{2,#1}}
\renewcommand{\qh}[4]{\overset{#4}{\mathsf{H}}{}^{#1,#2}_{#3}}
\renewcommand{\h}[3]{\overset{#3}{\mathsf{H}}{}^{#1}_{#2}}
\renewcommand{\pcoverset}[2]{\overset{#2}{\mathsf{C}}{}^{#1}}
\renewcommand{\qc}[3]{\overset{#3}{\mathsf{C}}{}^{#1,#2}}
\renewcommand{\srq}{\pr{-1/2,q}{}{}{}}
\renewcommand{\sdq}{\pdi{-1/2,q}{}{}{}}
\renewcommand{\gt}{\tau_{\tt t}}
\renewcommand{\gn}{\tau_{\tt n}}
\renewcommand{\chgt}{\check{\tau}_{\tt t}}
\renewcommand{\chgn}{\check{\tau}_{\tt n}}
\newcommand{\gtb}{\tau_{\tt t,b}}

\newcommand{\chgtb}{\check{\tau}_{\tt t,b}}

\renewcommand{\dom}{\Gamma}
\newcommand{\ML}{M_{\Lambda}}
\newcommand{\cML}{\calM_{\Lambda}}
\renewcommand{\omr}{\om_{\rho}}
\renewcommand{\Sr}{S_{\rho}}
\DeclareMathOperator{\divg}{div}
\renewcommand{\omb}{\Omega_{\tt b}}
\renewcommand{\domb}{\Gamma_{\tt b}}

\begin{document}

\maketitle{}

\begin{abstract}
We present both theory and an algorithm for solving time-harmonic wave problems in a general setting.
The time-harmonic solutions will be achieved by computing time-periodic solutions 
of the original wave equations. Thus, an exact controllability technique
is proposed to solve the time-dependent wave equations.
We discuss a first order Maxwell type system, which will be formulated 
in the framework of alternating differential forms. This enables us to investigate
different kinds of classical wave problems in one fell swoop, such as
acoustic, electro-magnetic or elastic wave problems. After a sufficient theory
is established, we formulate our exact controllability problem and
suggest a least-squares optimization procedure for its solution,
which itself is solved in a natural way by a conjugate gradient algorithm
operating in a purely $\lz$-type Hilbert space.
Therefore, it might be one of the biggest advances of this approach that 
the proposed conjugate gradient algorithm does not need any preconditioning.\\
\\
\amsclass{35Q60, 49M25, 65M99, 78A25, 78A30, 93B05, 93B40}
\end{abstract}

\newpage
\tableofcontents

\section{Introduction}\mylabel{introsec}

Time-harmonic wave propagation is an important phenomenon which has many obvious applications in acoustics,
electro-magnetics and elasticity, among others. Traditionally, the numerical solution approaches have
been based on finite differences, finite elements or boundary element techniques. As our goal is to
consider heterogeneous media as well, we pay attention to methods based on partial differential
equations. Hence, some kind of tessellation of the spatial domain is necessary.

To obtain accurate results for wave propagation, the discretization mesh needs to be adjusted to
the wavelength. If the time-harmonic case is directly addressed, one is faced with the solution
of a large-scale indefinite linear system which is a difficult task.

Instead of solving directly the time-harmonic problem for a given frequency $\omega\in\rzp$, 
it is possible to compute the solution by control techniques. 
Then the solution is found by searching for an appropriate initial data for the wave equation
which minimizes a quadratic functional that measures the difference between the initial state and the
final state after one time period $T=2\pi/\omega$. 
A natural quadratic error functional is the squared energy norm of the system,
allowing to minimize the cost by the conjugate gradient method (CGM) operating in Hilbert spaces.
This approach has been successfully applied to acoustics, electro-magnetics and elasticity
\cite{brigloe,brigloz,brigloro,glononlin,gloflow,glokellrein,glolioone,gloliotwo,gloro,heikeone,heiketwo,heiketre,heikefor}.
In practice, the method seems to have a good asymptotic computational
cost. Even though no theory exists, the computational cost of the method seems to be of order $\calO(n)$, 
where $n$ is the number of spatial degrees of freedom. 
The drawback of using the traditional (second order in time) formulation of the
wave equations is that the energy norm is then of $\h{1}{}{}$-type, and as such, the minimization problem is badly
conditioned. This is handled by applying preconditioning to the conjugate gradient minimization. Unfortunately,
this means that a discrete elliptic problem (linear system) still needs to be solved at every conjugate gradient
iteration step. In recent papers, the linear system has been solved by an algebraic multi-grid method which still
maintains the good asymptotic performance of the solution technique, but makes it quite more difficult to implement
the solver to utilize the computing power of modern parallel computers and multi-core processors.

Hence, an alternative approach has recently been proposed in the short paper of Glowinski et al. \cite{gloro}.
The idea is to formulate the control method for an equivalent first order system which has an $\lz$-type energy norm,
and hence, a well-conditioned minimization problem results. This eliminates the need for preconditioning the conjugate
gradient minimization and, thus, greatly simplifies the parallel implementation of the method. This approach also
has drawbacks as the spatial discretization needs to be based, for example, on mixed finite elements like Raviart-Thomas
elements, which are more difficult to implement than standard finite elements. Initial numerical experiments (still
unpublished) support the hypothesis that the cost of the new approach is also of order $\calO(n)$.

In our project, we aim at generalizing the approach of \cite{gloro} to generalized Maxwell equations formulated
in terms of differential forms. The same formulation covers electro-magnetic, acoustic and elastic waves and it can
be naturally discretized by so-called discrete differential forms (DDF) or discrete exterior calculus (DEC),
which has recently been under very active research \cite{hiranidiss,marsden,desmars}.
Our goal is to develop theory and software for efficiently solving
the generalized Maxwell equations using a control approach.
We present a new solution theory for the generalized Cauchy problem (CP) at hand, such that
we can be sure to have uniquely determined solutions evolving in time.
Here the papers \cite{weckcontrol,wecklip,paulytimeharm,paulystatic,paulydeco,paulyasym,kuhnpaulyreg}
as well as the monograph \cite{paulydiss} are useful.
Moreover, theoretical questions about the domain truncation procedure have to be considered.
We are planing to use absorbing boundary conditions (ABC), generalized Dirichlet-to-Neumann operators (DtN),
i.e. electric-to-magnetic operators (EtM), as well as perfectly matched layers (PML).
All these techniques have to be developed for differential forms.
The resulting software is targeted to mid-frequency variable
coefficient wave propagation problems in 2D and 3D domains, where the dimension of the computational domain is 10-100
wavelengths. The software is targeted for modern parallel computers and multi-core processors.

In this first paper we present and explain the basic ideas 
of our control approach for wave equations formulated as first order systems
and using generalized Maxwell equations in terms of differential forms.
First, in section \ref{cauchysec} we investigate the Cauchy problem and establish a solution theory which meets our needs 
utilizing the spectral calculus for unbounded selfadjoint linear operators in Hilbert space. 
Then, in section \ref{leastsquaressec} we introduce the least squares formulation 
and discuss the derivative of the least squares functional,
which is the essential ingredient in our resulting algorithm, since we plan to use a conjugate gradient method.
In section \ref{cgsec}, we discuss the conjugate gradient algorithm (CGA) in some detail.
In section \ref{classsec}, we translate our results presented in terms of differential forms to the classical framework of
vector analysis. We briefly demonstrate, which classical problems are covered by our generalized theory.
Finally, in section \ref{numresul} we present some preliminary numerical results 
and in section \ref{outlooksec} we outline the ongoing work in this project.

\section{Notations and preliminaries}\mylabel{notationsec}

We investigate wave scattering problems taking place in an exterior domain $\om$ of the Euclidean space $\rN$, 
which will be considered as a smooth $N$-dimensional differentiable Riemannian manifold 
with a Lipschitz boundary $\Gamma:=\p\Omega$.

We define the space $\cqunom$ of $\cu$-$q$-forms with compact support in $\om$.
This space admits a natural scalar product
$$\EH\mapsto\skp{E}{H}_\om:=\int_\om E\wedge*\bar{H}\in\cz,$$
where $*$ denotes the Hodge star operator with respect to the Euclidean metric in $\rN$, 
$\wedge$ the wedge product and the bar complex conjugation.
Using this scalar product and its induced norm we may define $\Lzqom{}$ as the closure of $\cqunom$.
Then $\Lzqom{}$ equipped with the scalar product
$\skp{\,\cdot\,}{\,\cdot\,}_{\Lzqom{}}:=\skp{\,\cdot\,}{\,\cdot\,}_\om$ becomes a Hilbert space,
the Hilbert space of square integrable $q$-forms on $\om$.

As usual, we denote the exterior derivative by $\rot$ 
and the co-derivative by $\pdiv$.
Thus we have on $q$-forms
$$\pdiv=(-1)^{(q-1)N}*\rot*.$$
With respect to the latter scalar product the linear operators $\rot$ and $\pdiv$
are formally skew-adjoint to each other,
i.e. for pairs of forms $\EH\in\cqunom\times\cqpeunom$ we have by the weak version of Stokes' theorem
\begin{align*}
0&=\intom\rot(E\wedge*\bar{H})=\intom\rot E\wedge*\bar{H}+(-1)^q\intom E\wedge\rot*\bar{H}\\
&=\intom\rot E\wedge*\bar{H}+\ub{(-1)^{q(N-q+1)}}_{=(-1)^{qN}}\intom E\wedge*\ub{*\rot*\bar{H}}_{=(-1)^{qN}\pdiv\bar{H}}\\
&=\skp{\rot E}{H}_{\lzqpeom}+\skp{E}{\pdiv H}_{\lzqom}.
\end{align*}
This yields the possibility for weak versions of $\rot$ and $\pdiv$ 
\big(in the sense of $\lzom$-valued distributions\big) 
using smooth, compactly supported forms as test-forms.
Hence, as usual we may define $\rot E$ for a $\lzqom$-form $E$ and say $E$ has weak exterior derivative, 
if there exists a $\lzqpeom$-form $G$, such that for all $\Phi\in\cqpeunom$ 
$$\skp{E}{\pdiv\Phi}_{\lzqom}=-\skp{G}{\Phi}_{\lzqpeom}$$
holds. Of course, we may define a weak co-derivative in the same way. Then we put
\begin{align*}
\rqom{}&:=\setb{E\in\Lzqom{}}{\rot E\in\Lzqpeom{}},\\
\dqom{}&:=\setb{H\in\Lzqom{}}{\pdiv H\in\qLzom{q-1}{}}.
\end{align*}
Equipped with their natural graph-norms these are Hilbert spaces.
Furthermore, we generalize the (electric) homogeneous boundary condition, 
which models a perfectly conducting obstacle and means 
that the tangential trace $\iota^*E$ of a differential form $E$ vanishes,
where $\iota:\dom\hookrightarrow\omq$ denotes the natural embedding of the boundary manifold $\dom$ 
regarded as an $(N-1)$-dimensional Riemannian submanifold of $\omq$.
For this purpose we define $\ronqom{}$ to be the closure of $\cqunom$ in the norm of $\rqom{}$.
Indeed, by Stokes' theorem and a density argument one may easily check
that for sufficiently smooth forms a vanishing tangential trace is generalized in $\ronqom{}$.
$\ronqom{}$ is also a Hilbert space as a closed subspace of $\rqom{}$.
An index $0$ at the lower left corners of the spaces $\ronqom{}$, $\rqom{}$
or $\dqom{}$ indicates vanishing exterior derivative or co-derivative, respectively.

Another way to define these Hilbert spaces is to look at the densely defined linear operator
$$\Abb{\ROT}{\cqunom\subset\lzqom}{\lzqpeom}{E}{\rot E}$$
and its adjoints, which will be marked by stars.
Then $\ol{\ROT}=\ROT^{**}$ is the weak exterior derivative on its domain of definition
$D(\ol{\ROT})=\ronqom{}$. The kernel of $\ol{\ROT}$ equals $\ronqnom{}$.
Its adjoint operator $\ROT^*=\ol{\ROT}^*$ equals by definition the negative weak co-derivative $-\DIV$
on its domain of definition $\dqpeom{}$, i.e.
$$\Abb{-\ROT^*=:\DIV}{\dqpeom{}\subset\lzqpeom}{\lzqom}{H}{\pdiv H}.$$
This is easy to see: Let $H\in D(\ROT^*)$ and $\ROT^*H=F$. Then by definition
$$\forall\,E\in D(\ROT)\quad\skp{\rot E}{H}_{\lzqpeom}=\skp{E}{F}_{\lzqom},$$
which is just the definition of the negative weak co-derivative. 
Therefore, $H$ is an element of $D(\ROT^*)=\dqpeom{}$ and $\ROT^*H=-\pdiv H$ holds.

Since $\pdiv\pdiv$ and $\rot\!\rot$ vanish in the smooth case,
$$\rot\!\rot=0,\quad\pdiv\pdiv=0$$
still hold true in the weak sense and we also have the well known and important formula
$$\rot\pdiv+\pdiv\rot=\Delta,$$
where the action of the Laplacian is to be understood componentwise with respect to Euclidean coordinates.
Moreover, we get with closures taken in $\lzqom$
$$\ol{\rot\pr{q-1}{}{}{(\circ)}(\om)}\subset\pr{q}{}{0}{(\circ)}(\om),\quad
\ol{\pdiv\pdi{q+1}{}{}{}(\om)}\subset\pdi{q}{}{0}{}(\om).$$

Let us formally define matrix-operators
$$M:=\zmat{0}{\pdiv}{\rot}{0},\quad\Lambda:=\zmat{\eps}{0}{0}{\mu},\quad\ML:=\ie\Lambda^\me M,$$
where $\eps$ respectively $\mu$ is a real, linear, symmetric, bounded and uniformly positive definite
\big(with respect to the $\Lzqom{}$- respectively $\Lzqpeom{}$-scalar product\big)
transformation on $q$- respectively $(q+1)$-forms, which is independent of time, and $\ie$ denotes the imaginary unit.
$\eps$ and $\mu$ model material properties,
i.e. in classical electro-magnetic theory $\eps$ is the dielectricity and $\mu$
the permeability of the underlying medium. We note that $\eps$ and $\mu$ are even allowed just to
have $\mathrm{L}^\infty(\om)$-entries in their matrix representations $\nu_{J',J}^h$ given by chart bases
$$\nu E=\sum_{J',J}\nu_{J',J}^hE_{J}\pd h_{J'}\quad\text{if }E=\sum_JE_{J}\pd h_{J}.$$

Since $\ol{\ROT}$ and $\DIV$ are skewadjoint to each other in this setting the unbounded linear operator
\beq\Abb{\cML}{D(\cML)\subset\Lzqqpeom{\Lambda}}{\Lzqqpeom{\Lambda}}{\EH}{{\ML}\EH}\mylabel{saOM}\eeq
with 
$$\Lzqqpeom{\Lambda}:=\Lzqqpeom{}:=\Lzqom{}\times\Lzqpeom{}$$
(as a set) equipped with the weighted scalar product
$$\skp{\,\cdot\,}{\,\cdot\,}_{\Lzqqpeom{\Lambda}}:=\skp{\Lambda\,\cdot\,}{\,\cdot\,}_{\Lzqqpeom{}},$$
where $\skp{\,\cdot\,}{\,\cdot\,}_{\Lzqqpeom{}}$ denotes the canonical scalar product in the product space $\Lzqqpeom{}$, 
and domain of definition
$$D(\cML):=\ronqom{}\times\dqpeom{}$$
is selfadjoint. The spectrum of $\cML$ might equal the entire real axis and we note
$${\ML}\EH=\ie(\eps^\me\pdiv H,\mu^\me\rot E).$$
For more details see \cite{paulydiss} or (for the classical case) \cite{leisbuch}.

\section{Problem formulation}\mylabel{formulationsec}

We are looking for $T$-periodic solutions in time of the following generalized Maxwell controllability problem
\begin{align}
(\p_t+\ie{\ML})\EH&=\FG&&\text{in }\Xi,\non\\
\gt E&=\lambda&&\text{in }\Upsilon,\mylabel{maxcontrol}\\
\EH(0)&=\EHnu&&\text{in }\om,\non\\
\EH(T)&\overset{!}{=}\EH(0)&&\text{in }\om,\non
\intertext{where $\gt$ denotes the tangential trace, i.e. $\gt=\iota^*$ in the smooth case.
Of course, the first equation may be written more explicitly as}
\p_tE-\eps^\me\pdiv H&=F&&\text{in }\Xi,\non\\
\p_tH-\mu^\me\rot E&=G&&\text{in }\Xi.\non
\end{align}
Here $I:=(0,T)$ with some time $T>0$ is an interval and $\bar{I}=[0,T]$ denotes its closure.
Furthermore, we introduce the two product sets $\Xi:=I\times\om$ and $\Upsilon:=I\times\dom$.

As a first order system, the problem at hand represents a natural generalization 
of classical wave equation problems associated to Helmholtz' equation.
In \cite{gloro} Glowinski et al. proposed an algorithm to compute time-$T$ periodic solutions $u$
of the prototypical scalar linear wave problem
\begin{align}
(\p_t^2-c^2\Delta)u&=0&&\text{in }\Xi,\non\\
\tau u&=g&&\text{on }\Upsilon,\mylabel{proto}\\
u(0)&=u(T)&&\text{in }\om,\non\\
\p_tu(0)&=\p_tu(T)&&\text{in }\om.\non
\end{align}
They utilized a truncation $\omr:=\setb{x\in\om}{|x|<\rho}$ 
of $\om$ introducing an artificial boundary 
(a sphere $\Sr$ of radius $\rho$ containing $\rN\ohne\om$)
and a first order absorbing boundary condition on it, i.e. setting the translation
of Sommerfeld's radiation condition to the time dependent formulation $(c^{\me}\p_{t}+\p_{r})u$ to zero.
Here $c$ is a positive real number and $g$ is a given time dependent boundary data. 
Furthermore, $\tau$ denotes the usual scalar trace operator and $r$ the Euclidean norm on $\rN$.

They transformed the latter system via the well known substitution
$$E:=\p_tu,\quad H:=\grad u$$
into a first order system of `linear acoustics'
\begin{align*}
\big(\p_t-\zmat{c^2}{0}{0}{1}\zmat{0}{\divg}{\grad}{0}\big)\EH&=(0,0)&&\text{in }I\times\omr,\\
\tau E&=\p_tg&&\text{on }\Upsilon,\\
c^{-1}E+\xi\cdot H&=0&&\text{on }I\times\Sr,\\
\EH(0)&=\EH(T)&&\text{in }\omr,
\end{align*}
which has a `Maxwell-type flavor', albeit simpler. Here $\xi(x):=x/|x|$.
One of the advantages of this first order system is that it allows for its solution an algorithm using
$$\lz(\omr)\times\lz(\omr)^N$$
as control space, i.e. the space of initial data. 
In former works there was always at least the first part of the control space a closed subspace of
$\h{1}{}{}(\omr)$, which makes the corresponding numerics much more difficult
due to the need of preconditioning, for instance, in conjugate gradient algorithms.
Such preconditioning is not necessary if one uses a purely $\lz(\omr)$-control space.

Utilizing the framework of alternating differential forms, our problem \eqref{maxcontrol} 
generalizes this approach not only to the classical Maxwell equations in three dimensions
but also to their generalized and coordinate free version.
We should mention that the generalized approach also comprises the system of linear acoustics 
and the 2-dimensional version of Maxwell's equations as well as the system of linear elasticity 
(with another boundary condition).

We emphasize that for $q:=0$ as well as $\FG:=(0,0)$, $\lambda:=g$, $\eps:=\mu:=1/c$ and $u:=E$
the original problem \eqref{proto} is recovered.

To start our analysis, we first have to establish a solution theory for the boundary value CP
\begin{align}
(\p_t+\ie{\ML})\EH&=\FG&&\text{in }\Xi,\non\\
\gt E&=\lambda&&\text{in }\Upsilon,\mylabel{cauchy}\\
\EH(0)&=\EHnu&&\text{in }\om\non
\end{align}
with given right hand sides $F$, $G$ and $\lambda$ as well as
initial data $\EHnu$ belonging to our control (Hilbert) space
$$\Hi:=\Lzqqpeom{\Lambda}.$$

\section{Solution theory for the Cauchy problem}\mylabel{cauchysec}

In order to solve \eqref{cauchy}, as a first step we must extend the boundary data from $\Upsilon$ to $\Xi$.

\subsection{Traces and extensions}

Recently Weck \cite{wecklip} showed how to obtain traces of differential forms on Lipschitz boundaries.
Let $\omb$ be a bounded Lipschitz domain in $\rN$ with boundary $\domb$. Then by \cite[Theorem 3]{wecklip}
there exists a linear and continuous tangential trace operator
$$\gtb:\prq{}(\omb)\To\srq(\domb),$$
where with the notations from \cite{wecklip}
$$\srq(\domb):=\setb{\lambda\in\qh{-1/2}{q}{\rho}{}(\domb)}{\Rotb\lambda\in\qh{-1/2}{q+1}{\rho}{}(\domb)}.$$
Here $\Rotb$ denotes the exterior derivative on $\domb$.
Moreover, by \cite[Theorem 4]{wecklip} $\gtb$ is surjective and there exists
a corresponding linear and continuous tangential extension operator (a right inverse of $\gtb$)
$$\chgtb:\srq(\domb)\To\prq{}(\omb).$$
Applying the well known Helmholtz decomposition
$$\Lzq{}(\omb)=\rot\pr{q-1}{}{}{\circ}(\omb)\oplus_\eps\dhqeps{}(\omb)\oplus_\eps\eps^\me\pdiv\dqpe{}(\omb),$$
where we introduce the finite dimensional space of Dirichlet forms
$$\dhqeps{}(\omb):=\ronqn{}(\omb)\cap\,\eps^\me\dqn{}(\omb),$$
and using $\gtb\ronq{}(\omb)=\{0\}$ we receive a linear and continuous tangential extension operator
$$\chgtb:\srq(\domb)\To\prq{}(\omb)\cap\,\eps^\me\pdiv\dqpe{}(\omb).$$
Now, we return to our exterior Lipschitz domain $\om\subset\rN$.
Using an usual cut-off-technique we obtain the following

\begin{lem}\mylabel{tracelem}
There exists a linear and continuous tangential trace operator
$$\gt:\rqom{}\To\srq(\dom)$$
and a corresponding linear and continuous tangential extension operator (right inverse)
$$\chgt:\srq(\dom)\To\rqom{}\cap\,\eps^\me\dqom{},$$
which even maps to forms with fixed \ul{com}p\ul{act} support and satisfies on $\srq(\dom)$
$$\gt\chgt=\id.$$
The kernel of $\gt$ equals $\ronqom{}$ and $\gt$ is even well defined on $\rqlocomq$. 
Moreover, $\chgt$ can be chosen, such that $\supp\chgt\lambda\subset\ol{\omr}$ 
holds for all $\lambda\in\srq(\dom)$ and for a fixed $\rho>0$ with $\rN\ohne\om\subset U_\rho$.
Here, $U_\rho\subset\rN$ denotes the open Euclidean ball with radius $\rho>0$ centered at the origin.
\end{lem}

\begin{rem}\mylabel{traceremh}
If the boundary is sufficiently smooth, e.g. $\pc{m+1}$, then even 
$$\gt E\in\qh{m-1/2}{q}{}{}(\dom)$$
holds for all forms $E\in\qhom{m}{q}{}{}$ or $E\in\qh{m}{q}{\loc}{}(\omq)$.
Moreover, applied to smooth forms from $\cqu(\omq)$ we have $\gt=\iota^*$ and,
of course, $\gt$ commutates 
with the exterior derivative. On the other hand, if $\lambda\in\qh{m-1/2}{q}{}{}(\dom)$
we may choose an extension, such that $\chgt\lambda\in\qhom{m}{q}{}{}$ holds and $\chgt\lambda$
is supported in $\ol{\omr}$. For details see \cite{kuhnpaulyreg}.
\end{rem}

$\gt$ and $\chgt$ may also be defined on time dependent forms. We get bounded linear operators
$$\gt:\mathsf{F}\big(I,\rqom{}\big)\To\mathsf{F}\big(I,\srq(\dom)\big)$$
and
$$\chgt:\mathsf{F}\big(I,\srq(\dom)\big)\To\mathsf{F}\big(I,\rqom{}\cap\,\eps^\me\dqom{}\big)$$
with similar properties as mentioned above, where the function space $\mathsf{F}$ could be, for instance,
$\pc{\ell}$, $\lz$, $\h{\ell}{}{}$.

Finally, we also need the corresponding normal trace and extension operators
$$\gn:=(-1)^{qN}\cast\gt*,\quad\chgn:=(-1)^{q(N-q)}*\chgt\cast$$
defined on $(q+1)$- respectively $(q-1)$-forms, where $\cast$ denotes Hodge's star operator on the
$(N-1)$-dimensional submanifold $\dom$ of $\omq$.

\subsection{Solution theory}

Now we return to the CP \eqref{cauchy}. Let $\lambda\in\pc{1}\big(\bar{I},\srq(\dom)\big)$. 
Then the canonical ansatz
$$\EHs:=\EH-(\chgt\lambda,0)$$
leads to a problem with homogeneous boundary condition
\begin{align}
(\p_t+\ie{\ML})\EHs&=\FGs&&\text{in }\Xi,\non\\
\gt\Es&=0&&\text{in }\Upsilon,\mylabel{homboundarycauchy}\\
\EHs(0)&=\EHsnu&&\text{in }\om\non
\end{align}
and new data
\begin{align*}
\FGs&:=\FG+(-\p_t\chgt\lambda,\mu^\me\rot\chgt\lambda),\\
\EHsnu&:=\EHnu-\big(\chgt\lambda(0),0\big).
\end{align*}
Since $\cML$ from \eqref{saOM} is linear and selfadjoint, 
spectral theory suggests a solution $\EHs$ of \eqref{homboundarycauchy} defined for all $t\in[0,\infty)$ by
\begin{align*}
\EHs(t)
&:=\exp(-\ie t\cML)\EHsnu+\int_0^t\exp\big(-\ie(t-s)\cML\big)\FGs(s)\,ds\\
&\;=\exp(-\ie t\cML)\Big(\EHsnu+\int_0^t\exp(\ie s\cML)\FGs(s)\,ds\Big).
\end{align*}
Let us analyze this solution thoroughly.
For instance, considering forms $\EHsnu\in\Hi$ and $\FGs\in\lz(I,\Hi)$ we obtain $\EHs$ in $\pc{0}(\bar{I},\Hi)$
and thus a solution
\beq\EH\in\pc{0}(\bar{I},\Hi),\mylabel{weaksol}\eeq
if $\EHnu\in\Hi$ and $\FG\in\lz(I,\Hi)$. Assuming stronger assumptions on the initial and right hand side data,
i.e. $\EHsnu\in D(\cML)$ and $\FGs\in\pc{0}(\bar{I},\Hi)\cap\lz\big(I,D(\cML)\big)$, 
we even get a solution $\EHs$ belonging to $\pc{1}(\bar{I},\Hi)\cap\pc{0}\big(\bar{I},D(\cML)\big)$.
Hence, we achieve a solution
$$\EH\in\pc{1}(\bar{I},\Hi)\cap\pc{0}\big(\bar{I},\rqom{}\times\dqpeom{}\big),$$
if, for instance,
\begin{align}
\EHnu&\in\rqom{}\times\dqpeom{},\non\\
\FG&\in\pc{0}(\bar{I},\Hi)\cap\lz\big(I,\rqom{}\times\dqpeom{}\big),\non\\
\rot\chgt\lambda(t)&\in\mu\dqpeom{},\quad t\in\bar{I},\mylabel{assstrongsol}\\
\gt F(t)&=\p_t\lambda(t),\non\\
\gt E_0&=\lambda(0).\non
\end{align}
Then $\EH$ is a solution of the CP \eqref{cauchy} in the strong sense.

Summing up we obtain:

\begin{theo}\mylabel{strongcauchytheo}
Let $\lambda\in\pc{1}\big(\bar{I},\srq(\dom)\big)$ as well as $\EHnu$ and $\FG$
satisfy \eqref{assstrongsol}. Then the CP \eqref{cauchy} is uniquely solved in 
$$\pc{1}(\bar{I},\Hi)\cap\pc{0}\big(\bar{I},\rqom{}\times\dqpeom{}\big)$$ 
by
\begin{align*}
\EH(t)&:=(\chgt\lambda,0)(t)+\exp(-\ie t\cML)\big(E_0-\chgt\lambda(0),H_0\big)\\
&\qquad+\int_0^t\exp\big(-\ie(t-s)\cML\big)(F-\p_s\chgt\lambda,G+\mu^\me\rot\chgt\lambda)(s)\,ds
\end{align*}
for $t\in\bar{I}$. We call $\EH$ the strong solution of the CP \eqref{cauchy} 
with data $(F,G,\lambda,E_0,H_0)$.
\end{theo}

Actually, we are interested in the purely $\lz$-type Hilbert space $\Hi$ as control space for the initial data 
and even not in $D(\cML)$ or $\rqom{}\times\dqpeom{}$.
Moreover, the constraints \eqref{assstrongsol} are too complicated and the assumptions on the data much too strong.
Thus, we have to weaken our solution concept. To approach weak solutions we first have to define suitable test forms.

\begin{defini}\mylabel{testforms}
For $\FPnu\in D(\cML)$ and $t\in\rz$ the family
$$\FP(t):=\exp(-\ie t\cML)\FPnu$$
defines a strong solution of the homogeneous Cauchy problem (HCP)
\begin{align*}
(\p_t+\ie{\ML})\FP&=(0,0)&&\text{in }\rz\times\om,\\
\gt\Phi&=0&&\text{in }\rz\times\dom,\\
\FP(0)&=\FPnu&&\text{in }\om.
\end{align*}
These solutions $\FP$ are elements of\, $\pc{1}(\rz,\Hi)\cap\pc{0}\big(\rz,D(\cML)\big)$
and we will call them test forms with initial values $\FPnu$.
\end{defini}

Next, we present the idea of the definition of weak solutions.
Thus, let $\EH$ be a strong solution of \eqref{cauchy} and
$\FP$ be a test form with initial value $\FPnu\in D(\cML)$. Then we may compute
\begin{align*}
\skpbH{\FG}{\FP}&=\skpbH{(\p_t+\ie{\ML})\EH}{\FP}\\
&=\skpbH{\p_t\EH}{\FP}-\skpb{M\EH}{\FP}_{\Lzqqpeom{}}\\
&=\p_t\skpbH{\EH}{\FP}-\skpbH{\EH}{\p_t\FP}\\
&\qquad-\skp{\rot E}{\Psi}_{\Lzqpeom{}}-\skp{\pdiv H}{\Phi}_{\Lzqom{}}.
\end{align*}
Since $\Phi\in\ronqom{}$ we obtain
$$\skp{\pdiv H}{\Phi}_{\Lzqom{}}=-\skp{H}{\rot\Phi}_{\Lzqpeom{}}$$
and assuming for these heuristic arguments that $E$, $\Psi$ and $\dom$ are sufficiently smooth we get by Stokes' theorem
\begin{align*}
&\qquad\skp{\rot E}{\Psi}_{\Lzqpeom{}}+\skp{E}{\pdiv\Psi}_{\Lzqom{}}
=\intom{\rot(E\wedge*\bar{\Psi})}\\
&=\int_\dom\iota^*(E\wedge*\bar{\Psi})
=(-1)^{qN}\int_\dom\iota^*E\wedge\cast\cast\iota^**\bar{\Psi}=\skp{\gt E}{\gn\Psi}_{\lzq(\dom)}.
\end{align*}
Putting all together yields
\begin{align*}
\skpbH{\FG}{\FP}&=\p_t\skpbH{\EH}{\FP}-\skp{\lambda}{\gn\Psi}_{\lzq(\dom)}\\
&\qquad-\skpbH{\EH}{\ub{(\p_t+\ie{\ML})\FP}_{=(0,0)}}.
\end{align*}
Hence, we only have to remove the time derivative from the forms $\EH$
to get our weak solutions. \big(Please compare to Weck \cite{weckcontrol}.\big)

\begin{defini}\mylabel{weakcauchydef}
Let $\EHnu\in\Hi$ and $\FG\in\lz(I,\Hi)$ as well as $\lambda\in\lz\big(I,\srq(\dom)\big)$.
Then the pair of forms $\EH$ is called a
weak solution of the CP \eqref{cauchy} with right hand side and initial data 
$(F,G,\lambda,E_0,H_0)$,
if and only if $\EH$ belongs to $\pc{0}(\bar{I},\Hi)$ and
\begin{align*}
&\qquad\skpbH{\EH}{\FP}(t)-\skpbH{\EHnu}{\FPnu}\\
&=\int_0^t\Big(\skpbH{\FG}{\FP}(s)+\skpdom{\lambda}{\gn\Psi}(s)\Big)\,ds
\end{align*}
holds for all $t\in\bar{I}$ as well as for all test forms $\FP$ with initial values $\FPnu\in D(\cML)$.
\end{defini}

\begin{rem}\mylabel{weakcauchyrem}
The term $\skpdom{\lambda}{\gn\Psi}$ needs some detailed interpretation.
The normal trace of a $(q+1)$-form from $\dqpeom{}$ is only an element of
$$\sdq(\dom):=\setb{\lambda\in\qh{-1/2}{q}{\pi}{}(\dom)}{\Div\lambda\in\qh{-1/2}{q-1}{\pi}{}(\dom)},$$
where $\Div:=(-1)^{(q-1)(N-1)}\circledast\Rot\circledast$ denotes the co-derivative on $\dom$
applied to $q$-forms. Please see again \cite{wecklip} for details.
Hence, at first sight the scalar product 
\beq\skpdom{\lambda}{\gn\Psi}(s)=\skpbdom{\lambda(s)}{\gn\Psi(s)}\mylabel{scalarprod}\eeq
for almost all $s$ makes only sense as an usual dual pairing
$$\gn\Psi(s)\lambda(s)=\skpb{\lambda(s)}{\gn\Psi(s)}_{\qh{1/2}{q}{\pi}{}(\dom),\qh{-1/2}{q}{\pi}{}(\dom)}.$$
Thus, $\lambda(s)$ should be an element of $\qh{1/2}{q}{\pi}{}(\dom)$ for almost all $s$, which is not the case.
But, since for almost all $s$ the boundary forms $\lambda(s)\in\srq(\dom)$ and $\gn\Psi(s)\in\sdq(\dom)$ 
have more `regularity' than $\qh{-1/2}{q}{\rho/\pi}{}(\dom)$, 
the scalar product \eqref{scalarprod} still makes sense for almost all $s$.
We will clarify this in the next lemma. 
\end{rem}

\begin{lem}\mylabel{lembilin}
The $\lzq(\dom)$-scalar product may be extended as a continuous bilinear form to 
$\srq(\dom)\times\sdq(\dom)$ (using Stokes' theorem) by the mapping
$$b:\srq(\dom)\times\sdq(\dom)\To\cz$$
with
$$b(\alpha,\beta)=\skp{\rot\chgt\alpha}{\chgn\beta}_{\lzqpeom}+\skp{\chgt\alpha}{\pdiv\chgn\beta}_{\lzqom}.$$
Moreover, for all $\EH\in\rqom{}\times\dqpeom{}$ Stokes' theorem 
$$\skp{\rot E}{H}_{\lzqpeom}+\skp{E}{\pdiv H}_{\lzqom}=b(\gt E,\gn H)$$
remains valid. Further on we will denote $b$ as usual by $\skp{\,\cdot\,}{\,\cdot\,}_{\lzq(\dom)}$.
\end{lem}

\begin{proof}
For $\alpha\in\srq(\dom)$ and $\beta\in\sdq(\dom)$ the respective extensions $\chgt\alpha$ and $\chgn\beta$ to $\om$
are elements of $\rqom{}$ and $\dqpeom{}$. Therefore, the definition of $b$ makes sense.
To show that $b$ is well defined, i.e. does not depend on the extensions, we pick some $\EH\in\rqom{}\times\dqpeom{}$
with $\gt E=\alpha$ and $\gn H=\beta$. Since $\gt(E-\chgt\alpha)$ and $\gn(H-\chgn\beta)$ vanish we have
$E-\chgt\alpha\in\ronqom{}$ and $H-\chgn\beta\in\donqpeom{}$. Thus, by definition (or an approximation argument) 
\begin{align*}
0&=\skpb{\rot(E-\chgt\alpha)}{H}_{\lzqpeom}+\skp{E-\chgt\alpha}{\pdiv H}_{\lzqom},\\
0&=\skp{\rot\chgt\alpha}{H-\chgn\beta}_{\lzqpeom}+\skpb{\chgt\alpha}{\pdiv(H-\chgn\beta)}_{\lzqom}
\end{align*}
hold. Addition shows $\skp{\rot E}{H}_{\lzqpeom}+\skp{E}{\pdiv H}_{\lzqom}=b(\alpha,\beta)$, 
which proves also the asserted formula.
Finally, the continuity of $b$ follows from the Cauchy-Scharz inequality and the continuity of the extensions.
\end{proof}

We are ready to prove the main result of this section.

\begin{theo}\mylabel{weakcauchytheo}
There exists at most one weak solution of \eqref{cauchy}. 
If additionally, for instance, $\lambda\in\h{1}{}{}\big(I,\srq(\dom)\big)$ 
then there exists always a unique weak solution of \eqref{cauchy},
which belongs to $\pc{0}\big([0,\infty),\Hi\big)$ since $T$ is arbitrary.
$[0,\infty)$ may be replaced by $\rz$ as well.
\end{theo}

\begin{proof}
The difference $\EH$ of two solutions satisfies 
$$\skpbH{\EH}{\FP}(t)=0$$
for all $t$ and all test forms $\FP$.
Since $\exp(\ie t\cML)$ is an unitary operator and $D(\cML)$ is dense in $\Hi$ we obtain 
$$\exp(\ie t\cML)\EH(t)=(0,0)$$ 
and thus $\EH(t)$ vanishes for all $t$, which proves uniqueness.
To show existence, we use the solution $\EH$ from Theorem \ref{strongcauchytheo}
suggested by spectral theory, which is still well defined and still belongs to $\pc{0}(\bar{I},\Hi)$ by \eqref{weaksol}
even with our weak assumptions. We note that we have replaced the stronger constraint 
$\lambda\in\pc{1}\big(\bar{I},\srq(\dom)\big)$ by the weaker constraint 
$\lambda\in\h{1}{}{}\big(I,\srq(\dom)\big)\subset\pc{0}\big(\bar{I},\srq(\dom)\big)$.
So, it remains to check if $\EH$ satisfies the integral equation of Definition \ref{weakcauchydef}.
For this purpose, let $\FP(t)=\exp(-\ie t\cML)\FPnu$, $t\in\rz$,
be a test form with $\FPnu\in D(\cML)$.
We start with the second term in the sum of the representation of $\EH$\,:
\begin{align*}
&\qquad\skpBH{\exp(-\ie t\cML)\big(E_0-\chgt\lambda(0),H_0\big)}{\FP(t)}\\
&=\skpBH{\big(E_0-\chgt\lambda(0),H_0\big)}{\FPnu}
=\skpbH{\EHnu}{\FPnu}-\skpb{\eps\chgt\lambda(0)}{\Phi_0}_{\Lzqom{}}
\end{align*}
The third term may be handled utilizing Fubini's theorem as follows:
\begin{align*}
&\qquad\skpBH{\int_0^t\exp\big(-\ie(t-s)\cML\big)\FGs(s)\,ds}{\FP(t)}\\
&=\int_0^t\skpbH{\exp(\ie s\cML)(F-\p_s\chgt\lambda,G+\mu^\me\rot\chgt\lambda)(s)}{\FPnu}\,ds\\
&=\int_0^t\skpbH{(F-\p_s\chgt\lambda,G+\mu^\me\rot\chgt\lambda)}{\FP}(s)\,ds\\
&=\int_0^t\skpbH{\FG}{\FP}(s)\,ds
+\int_0^t\skpbH{(-\p_s\chgt\lambda,\mu^\me\rot\chgt\lambda)}{\FP}(s)\,ds
\end{align*}
We proceed by calculating the last integral.
\begin{align*}
&\qquad-\int_0^t\skp{\p_s\chgt\lambda}{\eps\Phi}_{\Lzqom{}}(s)\,ds\\
&=-\int_0^t\p_s\skp{\chgt\lambda}{\eps\Phi}_{\Lzqom{}}(s)\,ds
+\int_0^t\skp{\chgt\lambda}{\eps\p_s\Phi}_{\Lzqom{}}(s)\,ds\\
&=-\skp{\chgt\lambda}{\eps\Phi}_{\Lzqom{}}(t)+\skpb{\chgt\lambda(0)}{\eps\Phi_0}_{\Lzqom{}}
+\int_0^t\skp{\chgt\lambda}{\pdiv\Psi}_{\Lzqom{}}(s)\,ds
\end{align*}
Hence, we get
\begin{align*}
&\qquad\int_0^t\skpbH{(-\p_s\chgt\lambda,\mu^\me\rot\chgt\lambda)}{\FP}(s)\,ds\\
&=-\skp{\chgt\lambda}{\eps\Phi}_{\Lzqom{}}(t)+\skpb{\chgt\lambda(0)}{\eps\Phi_0}_{\Lzqom{}}\\
&\qquad+\int_0^t\big(\ub{\skp{\chgt\lambda}{\pdiv\Psi}_{\Lzqom{}}(s)
+\skp{\rot\chgt\lambda}{\Psi}_{\Lzqpeom{}}(s)}_{\bds=\skpdom{\lambda}{\gn\Psi}(s)\eds}\big)\,ds
\end{align*}
by Lemma \ref{lembilin}. Putting all together completes the proof.
\end{proof}

\subsection{A new notation}

Let us change to a new and shorter notation,
which enables us to follow the forthcoming arguments and basic ideas more easily.
We set $\fn:=(0,0)$ as well as
\begin{align*}
\fu&:=\EH,&\ff&:=\FG,\\
\fun&:=\fu(0):=\EHnu,&\fel&:=(\chgt\lambda,0),\\
\fuT&:=\fu(T),&\fgl&:=-(\p_t+\ie \ML)\fel
=(-\p_t\chgt\lambda,\mu^\me\rot\chgt\lambda).
\end{align*}
With this notation our inhomogeneous Cauchy problem (ICP) \eqref{cauchy} reads as
\begin{align}
(\p_t+\ie{\ML})\fu&=\ff&&&&\text{in }\Xi,\non\\
\gt\pi\fu&=\lambda&&&&\text{in }\Upsilon,\mylabel{cauchyfu}\\
\fu(0)&=\fun&&&&\text{in }\om,\non
\intertext{where for a pair of forms $\pi$ denotes the projection onto the first component.
Moreover, $\fu$ may be decomposed into $\fu=\ful+\fuc$, where $\ful$ and $\fuc$
are the unique weak solutions of the CPs}
(\p_t+\ie{\ML})\ful&=\fn,&(\p_t+\ie{\ML})\fuc&=\ff&&\text{in }\Xi,\non\\
\gt\pi\ful&=0,&\gt\pi\fuc&=\lambda&&\text{in }\Upsilon,\mylabel{cauchyfulc}\\
\ful(0)&=\fun,&\fuc(0)&=\fn&&\text{in }\om.\non
\end{align}
$\ful$ depends linearly and continuously on the initial data $\fun$
and $\fuc$ is independent of the initial data $\fun$, i.e. constant with respect to $\fun$.
The unique weak solutions of \eqref{cauchyfu} and \eqref{cauchyfulc} exist by Theorem \ref{weakcauchytheo} 
in $\pc{0}\big(\bar{I},\Hi\big)$ for all $T$ and all
\beq\fun\in\Hi,\quad\ff\in\lz(I,\Hi),\quad\lambda\in\h{1}{}{}\big(I,\srq(\dom)\big)\mylabel{givendata}\eeq
and are given by the following formulas:
\begin{align}
\fu(t)&=\fel(t)+\e^{-\ie t\cML}\big(\fun-\fel(0)\big)+\int_0^t\e^{-\ie(t-s)\cML}(\ff+\fgl)(s)\,ds\non\\
\ful(t)&=\e^{-\ie t\cML}\fun\mylabel{specform}\\
\fuc(t)&=\fel(t)-\e^{-\ie t\cML}\fel(0)+\int_0^t\e^{-\ie(t-s)\cML}(\ff+\fgl)(s)\,ds\non
\end{align}

\section{Least-squares formulation of the controllability problem}\mylabel{leastsquaressec}

From now on, let the right hand side data $\ff$ and $\lambda$ satisfy \eqref{givendata}
as well as the time $T>0$ be given and fixed.

In order to solve the controllability problem \eqref{maxcontrol}, which reads now
\beq\text{\sf`Find }\fun\in\Hi\text{\sf, such that }\fu\text{\sf\, satisfies \eqref{cauchyfu} and }
\fuT=\fun\text{\sf.'}\mylabel{controlprob}\eeq
we investigate the equation
\beq\fuT-\fun=0\mylabel{eqfuTfun}\eeq
more thoroughly. With the help of \eqref{specform} we obtain
\begin{align*}
\fuT&=\fu(T)=\ful(T)+\fuc(T)=\e^{-\ie T\cML}\fun+\fucT,\\
\fuT-\fun&=(\e^{-\ie T\cML}-1)\fun+\fucT.
\end{align*}
Consequently, with the continuous linear operator in $\Hi$
$$\cC_t:=\cC(t):=\e^{-\ie t\cML}-1,$$
which satisfies $\norm{\cC_t}\leq2$ for all $t$ and will be called `control operator', we get
\beq\fuT-\fun=\cC_T\fun+\fucT.\mylabel{eqfuTfunC}\eeq
Hence, we have to solve the linear equation
$$\cC_T\fun+\fucT=0$$
in the Hilbert space $\Hi$, which we want to try approximately by an CGA.
Since, of course, $\cC_T$ is neither symmetric nor selfadjoint the usual CGA
suggests to consider the corresponding normal equation
\beq\cC_T^*\cC_T\fun+\cC_T^*\fucT=0,\mylabel{normaleq}\eeq
where $\cC_t^*=\e^{\ie t\cML}-1$ denotes the adjoint operator of $\cC_t$. We note $\cC_t^{**}=\cC_t$
and that $\cC_T^*\cC_T$ is selfadjoint. 
Consequently, we are forced to consider and to minimize the quadratic functional $\tilde{\cF}$ with
\begin{align*}
\tilde{\cF}(\fun)&:=\frac{1}{2}\skpH{\cC_T^*\cC_T\fun}{\fun}+\Re\skpH{\cC_T^*\fucT}{\fun}\\
&\,\,=\frac{1}{2}\skpH{\cC_T\fun}{\cC_T\fun}+\Re\skpH{\fucT}{\cC_T\fun}
=\frac{1}{2}\normH{\cC_T\fun+\fucT}^2-\frac{1}{2}\normH{\fucT}^2,
\end{align*}
which, of course, is minimized, if and only if the quadratic functional $\cF:=\tilde{\cF}+\normH{\fucT}^2/2$, i.e.
\beq\Abb{\cF}{\Hi}{[0,\infty)}{\fun}{\frac{1}{2}\normH{\cC_T\fun+\fucT}^2=\frac{1}{2}\normH{\fuT-\fun}^2},\mylabel{funcF}\eeq
is minimized. This leads to the following least-squares formulation:

{\sf`Find initial data $\fun\in\Hi$, such that}
\begin{align}
\forall\,\fvn\in\Hi\quad\cF(\fun)\leq\cF(\fvn).\text{\sf'}\mylabel{lsq}
\end{align}
Here, $\fu$ respectively $\fv$ is the unique weak solution of the ICP \eqref{cauchyfu}
with initial data $\fun$ respectively $\fvn$.

The implementation of the CGA in $\Hi$ is greatly facilitated by the knowledge of the derivative $\cF'$.
Since $\cF$ is differentiable as a quadratic functional we get from
\beq\cF(\fun+\fvn)=\cF(\fun)+\Re\skpH{\cC_T\fvn}{\cC_T\fun+\fucT}
+\frac{1}{2}\normH{\cC_T\fvn}^2,\mylabel{derical}\eeq
where $\fun,\fvn\in\Hi$, immediately
\begin{align}\begin{split}
\cF'(\fun)\fvn&=\Re\skpbH{\fvn}{\cC_T^*(\cC_T\fun+\fucT)}
=\Re\skpH{\fvn}{\cC_T^*\cC_T\fun+\cC_T^*\fucT}
\end{split}\mylabel{derivative}\end{align}
and, of course, the normal equation is recovered. 
In this sense, we may identify
$$\cF'(\fun)\qtext{with}\cC_T^*\cC_T\fun+\cC_T^*\fucT\in\Hi.$$
Furthermore, we receive the representations
\begin{align}
\cD_t&:=\cD(t):=\cC_t^*\cC_t
=(\e^{\ie t\cML}-1)(\e^{-\ie t\cML}-1)=2\big(1-\cos(t\cML)\big),\non\\
\fud_t&:=\fud(t):=\cC_t^*\fuc(t)\mylabel{cosrep}\\
&\,\,=(\e^{\ie t\cML}-1)\fel(t)+(\e^{-\ie t\cML}-1)\fel(0)
+\int_0^t(1-\e^{-\ie t\cML})\e^{\ie s\cML}(\ff+\fgl)(s)\,ds,\non
\end{align}
where we will call the continuous linear operator $\cD_t$ in $\Hi$ the `derivative operator'.
We have $\norm{\cD_t}\leq4$ for all $t$. Finally we obtain
\beq\cF'(\fun)\fvn=\Re\skpH{\fvn}{\cD_T\fun+\fud_T}.\mylabel{derivativetwo}\eeq
By \eqref{derical} we get also
$$\cF(\fun)\leq\cF(\fun+\fvn)-\cF'(\fun)\fvn$$
for all $\fvn\in\Hi$ and thus

\begin{rem}\mylabel{lsqequiFsn}
For $\fun\in\Hi$ the following assertions are equivalent:
\begin{itemize}
\item[\bf(i)]\quad $\fun$ is a solution of the least squares problem \eqref{lsq}.
\item[\bf(ii)]\quad $\cF'(\fun)=0$
\item[\bf(iii)]\quad $\cD_T\fun+\fud_T=0$ \big(normal equation \eqref{normaleq}\big)
\end{itemize}
\end{rem}

Using \eqref{eqfuTfunC}, let us interpret the derivative vector
$$\cD_T\fun+\fud_T=\cC_T^*\fusn\in\Hi,\quad\fusn:=\fuT-\fun\in\Hi$$
more thoroughly. Clearly, the forms $\fusp$ and $\fusm$ defined by
$$\fusp(t):=\e^{\ie t\cML}\fusn,\quad\fusm(t):=\e^{\ie(T-t)\cML}\fusn$$
are the unique weak solutions of the homogeneous adjoint Cauchy problems (HACPs$\pm$)
\begin{align}
(\p_t\mp\ie{\ML})\fuspm&=\fn&&&&\text{in }\Xi,\non\\
\gt\pi\fuspm&=0&&&&\text{in }\Upsilon,\mylabel{adjcauchy}\\
\fusp(0)&=\fusn,&\fusm(T)&=\fusn&&\text{in }\om\non
\end{align}
and we have $\fusp(T)=\fusm(0)=\e^{\ie T\cML}\fusn$, i.e.
$$\cC_T^*\fusn=(\e^{\ie T\cML}-1)\fusn=\fuspT-\fusn=\fusmn-\fusn.$$
Here, the signs $\pm$ indicate that the wave $\fusp$ evolves \ul{forward} in time,
whereas the wave $\fusm$ evolves \ul{backward} in time.
Of course, this implies a change of sign in the $\p_t$-term.
We note that we define the weak solutions of the HACPs$\pm$
analogously to Definition \ref{weakcauchydef}.
Finally, we obtain two more nice representations of our derivative vector utilizing 
the solutions of the HACPs$\pm$ \eqref{adjcauchy}
\beq\cD_T\fun+\fud_T=\fuspT-\fusn=\fusmn-\fusn.\mylabel{gradF}\eeq

As already pointed out, the derivative vector
depends on the initial condition $\fun$ both directly and indirectly
through the solution $\fu$ of the ICP \eqref{cauchyfu}
and one of the solutions $\fu^{*,\pm}$ of the HACPs$\pm$ \eqref{adjcauchy}.
Moreover, we saw in \eqref{specform} that $\fu=\ful+\fuc$ splits up into a linear and continuous
and a constant part (with respect to $\fun$).
Of course, the same holds true for the solutions of the adjoint equations.
Let us pick, for instance, the forward in time solution $\fus:=\fusp$.
Then $\fus$ depends linearly and continuously on the initial data $\fusn$
and may be decomposed into $\fus=\fusl+\fusc$, where $\fusl$ and $\fusc$
are the unique weak solutions of the HCPs
\begin{align*}
(\p_t-\ie{\ML})\fuslc&=\fn&&\text{in }\Xi,\\
\gt\pi\fuslc&=0&&\text{in }\Upsilon,\\
\fuslc(0)&=\fuslcn&&\text{in }\om
\end{align*}
with $\fusln:=\fulT-\fun=\cC_T\fun$ and $\fuscn:=\fucT$ as well as $\fusn=\fusln+\fuscn$.
Again, $\fusl$ depends linearly and continuously on $\fun$, whereas $\fusc$ does not depend on $\fun$.
Of course, we have
$$\fuslc(t)=\e^{\ie t\cML}\fuslcn.$$
Putting all together, we see
$$\cD_T\fun=\fuslT-\fulT+\fun,\quad\fud_T=\fuscT-\fucT.$$

\section{Conjugate gradient algorithm for the least-squares problem}\mylabel{cgsec}

Although it has become customary to use CGAs in Hilbert spaces, see e.g.
\cite{daniel,glononlin,gloflow,glokellrein,glolioone,gloliotwo,gloro} as a selection,
we briefly want to repeat the algorithm here.

In order to solve approximately our least squares problem (LSP) \eqref{lsq}, i.e. the linear equation
$$\cC_T\fun+\fucT=0$$
or by Remark \ref{lsqequiFsn} equivalently our normal equation \eqref{normaleq}
$$\cD_T\fun+\fud_T=\cC_T^*\cC_T\fun+\cC_T^*\fucT=0$$
we will use the following variant of the usual CGA:
Given an approximation $\fun^{n-1}$ and last search direction $\fd^{n-1}$ 
we compute the new search direction and approximation by
$$\fd^n:=\fr^{n-1}+\beta^{n-1}\fd^{n-1},\quad\fun^n:=\fun^{n-1}+\alpha^n\fd^n$$
with coefficients
\begin{align*}
\alpha^n&:=-\frac{\cF'(\fun^{n-1})\fd^n}{\normH{\cC_T\fd^n}^2}
=-\frac{\normH{\fr^{n-1}}^2}{\normH{\cC_T\fd^n}^2},\\
\beta^n&:=-\frac{\cF'(\fun^n)\cC_T^*\cC_T\fd^n}{\normH{\cC_T\fd^n}^2}
=-\frac{\normH{\fr^n}^2}{\alpha^n\normH{\cC_T\fd^n}^2}=\frac{\normH{\fr^n}^2}{\normH{\fr^{n-1}}^2},
\end{align*} 
where the residual is given by
$$\fr^n:=\cC_T^*\cC_T\fun^n+\cC_T^*\fucT=\fr^{n-1}+\alpha^n\cC_T^*\cC_T\fd^n.$$
We note that the initializing procedure in the CGA
\beq\fr^0=\cC_T^*\cC_T\fun^0+\cC_T^*\fucT,\mylabel{procone}\eeq	
i.e.
$$\fr^0=\cC_T^*\fusn=\fusT-\fusn,\quad\fusn=\cC_T\fun^0+\fucT=\fuT-\fun^0,$$
where we picked the forward in time solution $\fus:=\fu^{*,+}$,
needs the solution $\fu$ at time $T$ of the ICP \eqref{cauchyfu}
with initial data $\fun^0$ as well as the forward in time solution $\fus$ at time $T$
of the HACP+ \eqref{adjcauchy} with initial data $\fusn$. 
Analogously the procedure within the loop of the CGA 
\beq\fd=\cC_T^*\cC_T\fd^n,\mylabel{proctwo}\eeq
i.e.
$$\fd=\cC_T^*\fusn=\fusT-\fusn,\quad\fusn=\cC_T\fd^n=\fuT-\fd^n,$$
where we again used the forward in time solution $\fus:=\fu^{*,+}$,
needs the solution $\fu:=\ful$ at time $T$ of the HCP \eqref{cauchyfulc}
with initial data $\fd^n$ as well as the forward in time solution $\fus$ at time $T$
of the HACP+ \eqref{adjcauchy} with initial data $\fusn$.

We recall that the procedure \eqref{procone} respectively \eqref{proctwo}
may be identified with the calculation of the derivative or `gradient' 
$$\cF'(\fun^0)=\cF_{\fucT}'(\fun^0)\qtext{resp.}\cF_{0}'(\fd^n)$$
of the least squares functional 
$$\Abb{\cF:=\cF_{\fucT}}{\Hi}{[0,\infty)}{\fun}{\frac{1}{2}\normH{\cC_T\fun+\fucT}^2}$$
respectively
$$\Abb{\cF_{0}}{\Hi}{[0,\infty)}{\fun}{\frac{1}{2}\normH{\cC_T\fun}^2}.$$

We will present the CGA for the approximate solution of the LSP as our Algorithm \ref{algcgmain}.
In the beginning of the algorithm, before entering the iteration loop,
we choose an initial control vector $\fun^0\in\Hi$ and compute the first residual vector $\fr^0$, 
i.e. the `gradient' of the functional $\cF_{\fucT}$ at the point $\fun^0$,
which gives the first minimizing direction $\fd^1=\fr^0$.
The computation of this residual requires the solutions 
of the ICP \eqref{cauchyfu} with initial control vector $\fun^0$ 
and of the HACP+ \eqref{adjcauchy}.
Then, on each CGA iteration we calculate the solutions 
of the HCP \eqref{cauchyfulc} with initial vector $\fd^n$ 
and of the HACP+ \eqref{adjcauchy}.
This gives the `gradient' of the functional $\cF_{0}$ at the point $\fd^n$,
which is needed to update the new residual vector $\fr^n$ and the new control vector $\fun^n$.
Finally we set the new minimizing direction $\fd^{n+1}$.

\begin{algorithm}[H]
\caption{CGA in $\Hi$ for LSP \eqref{lsq}}
\mylabel{algcgmain}
\begin{algorithmic}
\STATE {\bf initialization}
\STATE {\bf set} $n=0$
\STATE {\bf set} initial control vector $\fun^n\in\Hi$
\STATE {\bf solve} ICP \eqref{cauchyfu} with initial vector $\fun^n$ and {\bf get} $\fu$
\STATE {\bf solve} HACP+ \eqref{adjcauchy} with initial vector $\fusn=\fuT-\fun^n$ and {\bf get} $\fus$
\STATE {\bf compute} residual vector \big(gradient $\cF_{\fucT}'(\fun^0)$\big) $\fr^n=\fusT-\fusn$
\STATE {\bf compute} norm $\rho^n=\normH{\fr^n}^2$
\IF    {$\rho^n$ small}
\STATE {\bf goto exit}
\ENDIF
\STATE {\bf set} first minimizing direction $\fd^{n+1}=\fr^n$
\LOOP  [for $n\geq1$ assuming $\fun^{n-1}$ and $\fr^{n-1}\neq0$, $\rho^{n-1}$ and $\fd^n\neq0$ are known]
\STATE {\bf solve} HCP \eqref{cauchyfulc} with initial vector $\fd^n$ and {\bf get} $\fu$
\STATE {\bf solve} HACP+ \eqref{adjcauchy} with initial vector $\fusn=\fuT-\fd^n$ and {\bf get} $\fus$
\STATE {\bf compute} gradient \big($\cF_{0}'(\fd^n)$\big) $\fd=\fusT-\fusn$
\STATE {\bf compute} parameter $\alpha=-\rho^{n-1}/\skpH{\fd}{\fd^n}$
\STATE {\bf update} control vector $\fun^n=\fun^{n-1}+\alpha\fd^n$
\STATE {\bf update} residual vector $\fr^n=\fr^{n-1}+\alpha\fd$
\STATE {\bf compute} norm $\rho^n=\normH{\fr^n}^2$
\IF    {$\rho^n$ small {\bf or} $n$ large}
\STATE {\bf goto exit}
\ENDIF
\STATE {\bf compute} parameter $\rho=1/\rho^{n-1}$
\STATE {\bf compute} parameter $\rho=\rho^n\rho$
\STATE {\bf update} minimizing direction $\fd^{n+1}=\fr^n+\rho\fd^n$
\STATE {\bf set} $n=n+1$
\ENDLOOP
\STATE {\bf exit}
\STATE {\bf take} $\fun^n$ as solution
\end{algorithmic}
\end{algorithm}

We note that we may use the backward in time system HACP- \eqref{adjcauchy}
instead of HACP+ as well. Then, in this variant by \eqref{gradF} 
we have to replace the computation of the residual or gradient vector 
$\fusT-\fusn=\fuT^{*,+}-\fusn$ by $\fun^{*,-}-\fusn$.

\section{Translation to classical problems}\mylabel{classsec}

We briefly mention, which classical problems of vector analysis are covered by our general CP \eqref{cauchy}.
\eqref{cauchy} reads:
\begin{align}
\p_tE-\eps^\me\pdiv H&=F&&\text{in }\Xi\non\\
\p_tH-\mu^\me\rot E&=G&&\text{in }\Xi\mylabel{cauchydetail}\\
\gt E&=\lambda&&\text{in }\Upsilon\non
\end{align}
The initial condition always stays the same.

In $\rz^3$ the exterior derivative $\rot$ and the co-derivative $\pdiv$ 
turn to the classical differential operators from vector analysis
$$\grad=\nabla,\quad\curl=\nabla\,\times\,,\quad\divg=\nabla\,\cdot\,$$
and the well known standard Sobolev spaces 
$$\lzom,\quad\h{}{}{(\circ)}(\grad,\om)=\h{1}{}{(\circ)}(\om),\quad
\h{}{}{(\circ)}(\curl,\om),\quad\h{}{}{(\circ)}(\divg,\om)$$
appear. Moreover, the tangential trace becomes the usual scalar, tangential or normal trace, respectively.
As long as the operators $\times$ or $\curl$ are not involved, the classical calculus extends to $\rN$, $N\in\nz$. 

We obtain the following problems in $\rN$\,:\\
$\bullet$ $q=N$ (trivial case): $\p_tE=F$\\
$\bullet$ $q=0$ (linear acoustics, Dirichlet case):
\begin{align*}
\p_tE-\eps^\me\divg H&=F&&\text{in }\Xi\\
\p_tH-\mu^\me\grad E&=G&&\text{in }\Xi\\
\restr{E}{\dom}&=\lambda&&\text{on }\Upsilon
\intertext{$\bullet$ $q=N-1$  (linear acoustics, Neumann case):}
\p_tE-\eps^\me\grad H&=F&&\text{in }\Xi\\
\p_tH-\mu^\me\divg E&=G&&\text{in }\Xi\\
\nu\cdot\restr{E}{\dom}&=\lambda&&\text{on }\Upsilon
\intertext{$\bullet$ $q=1$ and $N=3$  (Maxwell's equations):}
\p_tE+\eps^\me\curl H&=F&&\text{in }\Xi\\
\p_tH-\mu^\me\curl E&=G&&\text{in }\Xi\\
\nu\times\restr{E}{\dom}&=\lambda&&\text{on }\Upsilon
\end{align*}

We note that the equations of linear elasticity are also covered by our approach,
if we change the tangential boundary condition into the more simple one
of componentwise scalar Dirichlet boundary conditions. 
See, for example, Weck and Witsch \cite{linelae}.

\section{Conclusion and outlook}\mylabel{outlooksec}

The considered approach of combining the exact controllability method with DEC-based
discretization seems to be a promising way to compute time-harmonic scattered waves. 
The problem setup is general enough to treat the most important cases of linear wave propagation, 
i.e. electro-magnetic, acoustic and elastic waves in three space dimensions.

There are certain key benefits of the method. The DEC approach leads to a discrete scheme that
has good conservation properties \cite{marsden}. Also, the resulting time integration scheme
can be implemented in explicit manner. 
As the mass matrices are also diagonal, the time integration will be
computationally very efficient and all the related computations are easily parallelized by
standard domain decomposition techniques with coarse grained boundary swapping approach.
The parallel implementation of the outer CG-iterations is also very straightforward because
we minimize the error in periodicity in the squared energy norm of the system. In the family
of problems we consider, the energy norm is a weighted $\lz$-norm, which means that the
discrete quadratic functional we minimize is spanned by a diagonal mass matrix. Hence, in
practice, no preconditioning is required. This is supported by our initial experiments (appendix), 
at least when the mesh is refined. The convergence of CG iterations did not depend on the mesh step size.

The use of DEC causes some challenges to overcome. These are related to the definition of
the dual mesh and the resulting discrete Hodge operator. In our initial experiments (appendix), we used
well-centered meshes and the circumcenter based dual mesh definition, which naturally
leads to a diagonal discrete Hodge operator. The tiling of a general domain with well-centered
simplices is an open problem, which has been accomplished for some simple shaped domains \cite{vanderzee}.

To make the method applicable to general cases arising from practical problems, we must allow for
simplicial meshes, which are not well-centered. Therefore, e.g., the barycentric dual meshes
need to be considered. It should be noted that simplicial meshes are not obligatory in our approach.
As was shown in \cite{marsden}, the classical and very widely used Yee-scheme for Maxwell's equations 
\cite{yee} is just a special case of general DEC-approach and the same control scheme can be implemented
directly for Yee's scheme as well.

\appendix

\section{Appendix: Preliminary numerical results}\mylabel{numresul}

We have implemented the DEC for 3-dimensional geometries to solve electro-magnetic problems. 
Following \cite[section 4.3]{marsden} the implementation allows us 
to have an unstructured mesh in space and asynchronous time steps. 
Our implementation is based on the circumcentric dual mesh, 
which is the most simple way to build a DEC solver, 
but it also requires Delaunay's property of the mesh. 
Since we are interested in time-periodic solutions, 
we implemented the CGA using the theory of section \ref{cgsec} as well.

\begin{figure}[htp]
\centering
\includegraphics{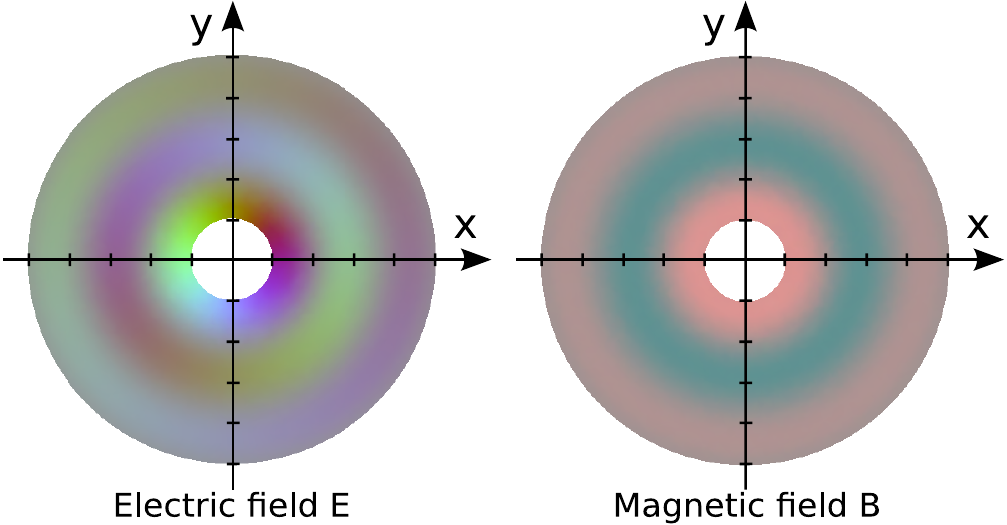}
\caption{We see simulated fields at $t=T$. The pictures are cross sections on a $xy$-plane, 
where the $x$-, $y$- and $z$-components of the field vectors are presented 
in blue, green and red, respectively. 
The zero field would be $50\%$ gray.}
\mylabel{fig_balls}
\end{figure}

We discuss some preliminary results of our simulations. 
Let us consider a scattering problem, where electro-magnetic plane waves 
are reflected by a sphere and scattered to infinity. 
We are interested in the accuracy of the simulation
and in the convergence of the CGA.

\begin{figure}[htp]
\centering
\includegraphics{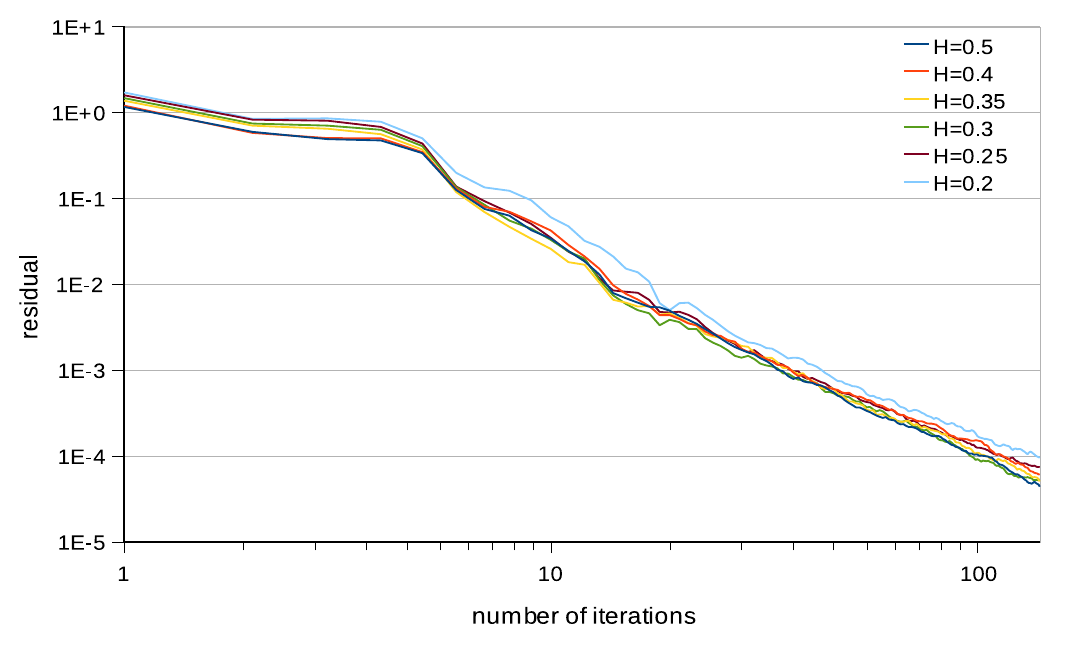}
\caption{convergence of the residual in CGA}
\label{fig_convergence}
\end{figure}
\begin{figure}[htp]
\centering
\includegraphics{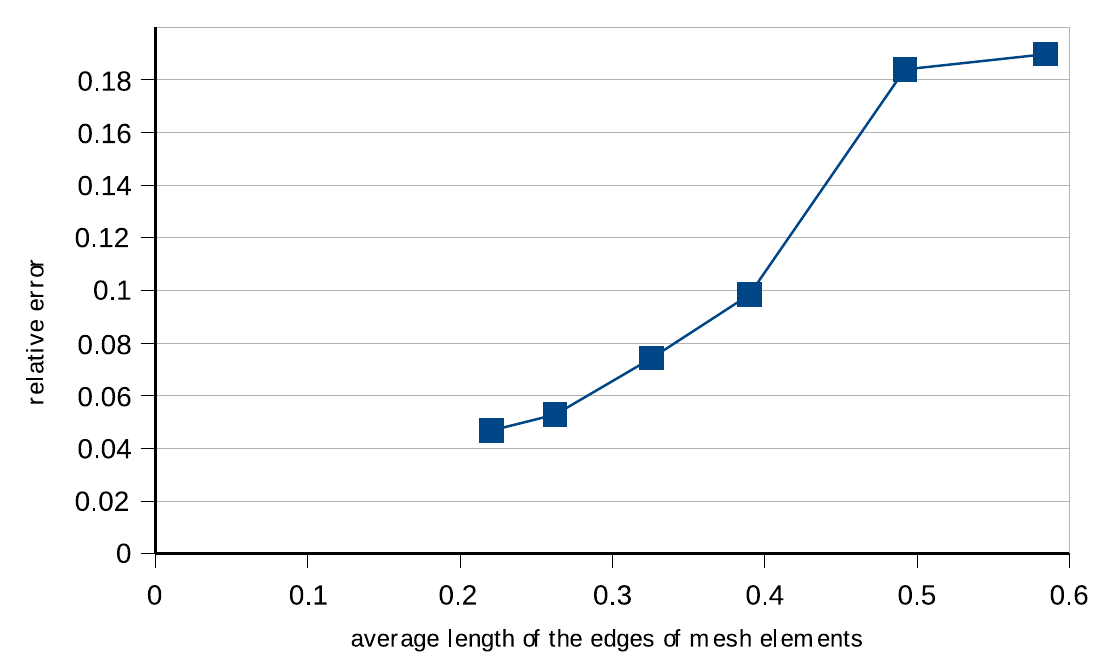}
\caption{relative error of the simulated electric field at time $t=T$ integrated over the mesh volume}
\label{fig_error}
\end{figure}

For our very simple model radiation problem 
\begin{align*}
\curl H-\ie\omega E&=0&&\text{in }\Omega,\\
\curl E+\ie\omega H&=0&&\text{in }\Omega,\\
\nu\times\restr{E}{\dom}&=\lambda&&\text{on }\Gamma,
\end{align*}
where $\Omega:=\setb{x\in\rd}{|x|>1}$ is the exterior of the closed unit ball, 
i.e., $\rd$ with a spherical hole (ball) of radius $1$ in the middle,
$\Gamma=S^2$ and where we picked the frequency $\omega=2\pi/3$,
the exact solution is known explicitly and can be found, for instance,
in the book \cite[Theorem 6.25]{colton}. 
We took only one non-zero component $a_1^0=1$. 
This simplifies the solution to
$$H=\frac{\ie}{\omega}\curl E,\quad E=\curl\Es,\quad\Es=h(\omega r)y(\xi)\id$$
with Hankel's spherical function of first kind $h:=h_1^{(1)}$ and spherical harmonic $y:=y_1^0$.
More explicitly, we have 
$$E=h(\omega r)Y(\xi)\times\xi,\quad Y:=\tilde{\nabla}y,$$ 
where $\tilde{\nabla}$ denotes the spherical gradient on $\Gamma$.
Let us note that
$$\nu\times E=h(\omega r)\xi\times\xi\times Y(\xi)=h(\omega r)Y(\xi)$$
since $\nu=-\xi$ and $Y(\xi)$ is tangential at $\Gamma$.
We generate a wave on the boundary $\Gamma$ by the Dirichlet boundary condition
$\lambda:=h(\omega)Y$, i.e., setting $r:=1$.
Picking an artificial outer boundary $\tilde{\Gamma}$, the sphere of radius $5$ centered at the origin, 
we impose the classical Silver-M\"uller first order absorbing boundary condition
$$\xi\times E+H=0\quad\text{on }\tilde{\Gamma}.$$ 

We have simulated the test problem with six different meshes of varying element sizes. 
The initialized edge lengths varied from about $1/5$ to $1/2$. 
In Figure \ref{fig_convergence} we see how the residual converges
in the CGA. After $140$ loops of the CGA
we got about $10^{-5}$ times smaller residuals.
The convergence seems to be independent of the mesh element size. 
In Figure \ref{fig_error} we plotted the differences between the simulation 
corresponding to different mesh element sizes and the exact solution. 
The error of the simulated fields is decreasing when the mesh is refined. 
The decrease of the error even might be of second order 
with respect to the average edge length.

\begin{acknow}
The first author expresses his gratitude to the Department of Mathematical Information Technology
of the University of Jyv\"askyl\"a (Finland) for scientific and financial support.
The authors thank Jukka R\"abin\"a (Jyv\"askyl\"a) 
for providing the results and pictures of the appendix. 
\end{acknow}

\vspace*{2cm}
\begin{center}
\begin{tabular}{ll}
{\sf Dirk Pauly} & {\sf Tuomo Rossi}\\
\\ 
Universit\"at Duisburg-Essen & University of Jyv\"askyl\"a\\ 
Campus Essen &\\
Fakult\"at f\"ur Mathematik & Faculty of Information Technology\\
 & Department of\\
 & Mathematical Information Technology\\
Universit\"atsstr. 2 & P.O. Box 35 (Agora)\\
45117 Essen & FI-40014 Jyv\"askyl\"a\\
Germany & Finland\\ 
e-mail: dirk.pauly@uni-due.de & e-mail: tuomo.rossi@jyu.fi
\end{tabular}
\end{center}

\end{document}

%% file: head_paper.tex
\usepackage{a4,exscale,ifthen,amsfonts,amssymb,amsmath,amscd,graphicx}
\usepackage[english]{babel}
\usepackage[mathscr]{eucal}

\author{{\sf\pauthor}}
\numberwithin{equation}{section}

\newenvironment{acknow}{{\vspace*{1cm}\noindent\bf Acknowledgements }}{}
\newcommand{\bewboxw}{\mbox{}\hfill $\square$ \\}

\newenvironment{proof}{{\noindent\bf Proof }}{\bewboxw}

\newcommand{\keywords}[1]{{\noindent\bf Key Words }#1}
\newcommand{\amsclass}[1]{{\noindent\bf AMS MSC-Classifications }#1}

\ifthenelse{\equal{\mylabelonoff}{on}}
{\newcommand{\mylabel}[1]{\label{#1}\fbox{{\rm #1}}}}{\newcommand{\mylabel}[1]{\label{#1}\makebox[0mm][]{}}}
\ifthenelse{\equal{\allowdisbrk}{yes}}
{\allowdisplaybreaks}{}

\newcommand{\paper}[7]{\bibitem{#1} #2, `#3', {\it #4}, #5, (#6), #7.}
\newcommand{\paperpub}[5]{\bibitem{#1} #2, `#3', {\it #4}, (#5), accepted for publication.}

\newcommand{\papersub}[4]{\bibitem{#1} #2, `#3', (#4), submitted.}
\newcommand{\book}[6]{\bibitem{#1} #2, {\it #3}, #4, #5, (#6).}

\newcommand{\dissavail}[6]{\bibitem{#1} #2, `#3', {\sf Dissertation}, #4, (#5), available from {\tt #6}.}

\newcommand{\repjyu}[8]{\bibitem{#1} #2, `#3', {\it Reports of the Department of Mathematical Information Technology}, 
University of Jyv\"askyl\"a, Series #4. Scientific Computing, No. #4. #5/#6, ISBN #7, ISSN #8.}
\newcommand{\reportavail}[7]{\bibitem{#1} #2, `#3', {\it #4}, #5, (#6), available from {\tt #7}.}
\newcommand{\preprintarxiv}[5]{\bibitem{#1} #2, `#3', (#4), preprint, available at {\tt arXiv}: #5.}

%% file: pauly_kopf.tex
\usepackage{amsfonts,amssymb,amsmath,amscd}
\usepackage[mathscr]{eucal}

\newcommand{\schluss}{\ifodd\value{page}\newpage\thispagestyle{empty}\makebox[0mm][]{}\color{sehrhell}.\fi\end{document}}
\newcommand{\non}{\nonumber}
\newcommand{\normabst}{\hspace{-0.4ex}}

\newcommand{\ol}[1]{\overline{#1}}
\newcommand{\ul}[1]{\underline{#1}}

\newcommand{\ub}{\underbrace}
\newcommand{\qtext}[1]{\quad\text{#1}\quad}

\newcommand{\cast}{\circledast}

\newcommand{\me}{{-1}}

\newcommand{\intom}{\int_\Omega}

\newcommand{\om}{{\Omega}}

\newcommand{\omq}{\ol{\om}}
\newcommand{\omb}{\om_\mathrm{b}}

\newcommand{\dom}{{\p\om}}
\newcommand{\domb}{{\p\omb}}

\newcommand{\ohne}{\setminus}

\newcommand{\eps}{\varepsilon}

\newcommand{\EH}{(E,H)}

\newcommand{\EHs}{(\tilde{E},\tilde{H})}

\newcommand{\FG}{(F,G)}

\newcommand{\FGs}{(\tilde{F},\tilde{G})}

\newcommand{\To}{\longrightarrow}

\newcommand{\restr}[2]{\left.#1 \right|_{#2}}

\newcommand{\Abb}[5]{\begin{array}{ccccc}#1&:&#2&\longrightarrow&#3\\{}&{}&#4&\longmapsto&#5\end{array}}

\newcommand{\bds}{\begin{displaystyle}}
\newcommand{\eds}{\end{displaystyle}}
\newcommand{\beq}{\begin{equation}}
\newcommand{\eeq}{\end{equation}}
\newcommand{\beqa}{\begin{eqnarray}}
\newcommand{\eeqa}{\end{eqnarray}}
\newcommand{\beqanon}{\begin{eqnarray*}}
\newcommand{\eeqanon}{\end{eqnarray*}}
\newcommand{\bamath}{$$\begin{array}}
\newcommand{\eamath}{\end{array}$$}
\newcommand{\ba}{\begin{array}}
\newcommand{\ea}{\end{array}}
\newcommand{\bc}{\begin{center}}
\newcommand{\ec}{\end{center}}
\newcommand{\bmini}{\begin{minipage}}
\newcommand{\emini}{\end{minipage}}

\newtheorem{lem}{Lemma}[section]
\newtheorem{theo}[lem]{Theorem}

\newtheorem{defini}[lem]{Definition}
\newtheorem{rem}[lem]{Remark}

\newcommand{\setb}[2]{\big\{#1 \;\text{\bf :}\;  #2\big\}}

\newcommand{\rz}{\mathbb{R}}

\newcommand{\nz}{\mathbb{N}}

\newcommand{\cz}{\mathbb{C}}

\newcommand{\rd}{{\rz^3}}
\newcommand{\rN}{{\rz^N}}

\newcommand{\rzp}{{\rz_+}}

\newcommand{\calO}{{\mathscr{O}}}

\newcommand{\calH}{{\mathscr{H}}}

\newcommand{\calM}{{\mathscr{M}}}

\DeclareMathOperator{\cF}{{\mathcal{F}}}

\newcommand{\zmat}[4]{\begin{bmatrix}#1&#2\\#3&#4\end{bmatrix}}

\DeclareMathOperator{\p}{\partial}

\DeclareMathOperator{\loc}{loc}

\DeclareMathOperator{\id}{Id}

\DeclareMathOperator{\supp}{supp}

\DeclareMathOperator{\rot}{rot}
\DeclareMathOperator{\Rot}{Rot}
\DeclareMathOperator{\ROT}{ROT}
\DeclareMathOperator{\pdiv}{div}
\DeclareMathOperator{\Div}{Div}
\DeclareMathOperator{\DIV}{DIV}

\DeclareMathOperator{\ie}{i}
\DeclareMathOperator{\e}{e}

\DeclareMathOperator{\pd}{d\hspace{-0.4ex}}

\DeclareMathOperator{\curl}{curl}

\newcommand{\pcoverset}[2]{\overset{#2}{\mathrm{C}}{}^{#1}}
\newcommand{\pc}[1]{\pcoverset{#1}{}}

\newcommand{\cu}{\pc{\infty}}

\newcommand{\lz}{\mathrm{L}^2}
\newcommand{\lzom}{\lz(\Omega)}

\newcommand{\h}[3]{\overset{#3}{\mathrm{H}}{}^{#1}_{#2}}

\newcommand{\qc}[3]{\overset{#3}{\mathrm{C}}{}^{#1,#2}}
\newcommand{\cq}[1]{\qc{#1}{q}{}}
\newcommand{\cqn}[1]{\qc{#1}{q}{\circ}}
\newcommand{\cqu}{\cq{\infty}}
\newcommand{\cqun}{\cqn{\infty}}

\newcommand{\cqunom}{\cqun(\Omega)}

\newcommand{\cqpen}[1]{\qc{#1}{q+1}{\circ}}

\newcommand{\cqpeun}{\cqpen{\infty}}

\newcommand{\cqpeunom}{\cqpeun(\Omega)}

\newcommand{\qlz}[1]{\mathrm{L}^{2,#1}}
\newcommand{\qlzom}[1]{\qlz{#1}(\Omega)}
\newcommand{\lzq}{\qlz{q}}
\newcommand{\lzqom}{\qlzom{q}}
\newcommand{\qLz}[2]{\qlz{#1}_{#2}}
\newcommand{\qLzom}[2]{\qLz{#1}{#2}(\Omega)}
\newcommand{\Lzq}[1]{\qLz{q}{#1}}
\newcommand{\Lzqom}[1]{\Lzq{#1}(\Omega)}

\newcommand{\lzqpeom}{\qlzom{q+1}}
\newcommand{\Lzqpe}[1]{\qLz{q+1}{#1}}
\newcommand{\Lzqpeom}[1]{\Lzqpe{#1}(\Omega)}

\newcommand{\qh}[4]{\overset{#4}{\mathrm{H}}{}^{#1,#2}_{#3}}

\newcommand{\qhom}[4]{\qh{#1}{#2}{#3}{#4}(\Omega)}

\newcommand{\pr}[4]{{}_{#3}\overset{#4}{\mathrm{R}}{}^{#1}_{#2}}
\newcommand{\prq}[1]{\pr{q}{#1}{}{}}

\newcommand{\rqloc}{\pr{q}{\loc}{}{}}

\newcommand{\ronq}[1]{\pr{q}{#1}{}{\circ}}

\newcommand{\ronqn}[1]{\pr{q}{#1}{0}{\circ}}

\newcommand{\rqom}[1]{\prq{#1}(\Omega)}

\newcommand{\rqlocomq}{\rqloc(\ol{\Omega})}

\newcommand{\ronqom}[1]{\ronq{#1}(\Omega)}

\newcommand{\ronqnom}[1]{\ronqn{#1}(\Omega)}

\newcommand{\pdi}[4]{{}_{#3}\overset{#4}{\mathrm{D}}{}^{#1}_{#2}}
\newcommand{\pdq}[1]{\pdi{q}{#1}{}{}}
\newcommand{\dqpe}[1]{\pdi{q+1}{#1}{}{}}

\newcommand{\dqn}[1]{\pdi{q}{#1}{0}{}}

\newcommand{\donqpe}[1]{\pdi{q+1}{#1}{}{\circ}}

\newcommand{\dqom}[1]{\pdq{#1}(\Omega)}
\newcommand{\dqpeom}[1]{\dqpe{#1}(\Omega)}

\newcommand{\donqpeom}[1]{\donqpe{#1}(\Omega)}

\newcommand{\dH}[3]{{}_{#3}\calH^{#1}_{#2}}

\newcommand{\skp}[2]{\langle#1,#2\rangle}

\newcommand{\skpb}[2]{\big\langle#1,#2\big\rangle}
\newcommand{\skpB}[2]{\Big\langle#1,#2\Big\rangle}

\newcommand{\norm}[1]{|\normabst|#1|\normabst|}

\newcommand{\qqLzom}[3]{\qLzom{#1,#2}{#3}}
\newcommand{\Lzqqpeom}[1]{\qqLzom{q}{q+1}{#1}}

\DeclareMathOperator{\grad}{grad}

\newcommand{\dhqeps}{{\dH{q}{}{\eps}}}

\usepackage{algorithm,algorithmic}

\newcommand{\Hi}{\mathbb{H}}

\newcommand{\Es}{\tilde{E}}

\newcommand{\EHnu}{(E_0,H_0)}

\newcommand{\EHsnu}{(\tilde{E}_0,\tilde{H}_0)}
\newcommand{\FP}{(\Phi,\Psi)}
\newcommand{\FPnu}{(\Phi_0,\Psi_0)}

\newcommand{\skpH}[2]{\skp{#1}{#2}_{\Hi}}
\newcommand{\skpbH}[2]{\skpb{#1}{#2}_{\Hi}}
\newcommand{\skpBH}[2]{\skpB{#1}{#2}_{\Hi}}
\newcommand{\skpdom}[2]{\skp{#1}{#2}_{\Lzq{}(\dom)}}
\newcommand{\skpbdom}[2]{\skpb{#1}{#2}_{\Lzq{}(\dom)}}

\newcommand{\normH}[1]{\norm{#1}_{\Hi}}

\newcommand{\gt}{\gamma_\tau}
\newcommand{\gn}{\gamma_\nu}
\newcommand{\chgt}{\check{\gamma}_\tau}
\newcommand{\chgn}{\check{\gamma}_\nu}
\newcommand{\srq}{\cR^q}
\newcommand{\sdq}{\mathcal{D}^q}

\newcommand{\fu}{\mathbf{u}}
\newcommand{\fun}{\fu_0}
\newcommand{\fus}{\fu^*}
\newcommand{\fusn}{\fus_0}
\newcommand{\fusT}{\fus_T}
\newcommand{\ff}{\mathbf{f}}
\newcommand{\pfg}{\mathbf{g}}

\newcommand{\fe}{\mathbf{e}}
\newcommand{\fel}{\fe_\lambda}
\newcommand{\fgl}{\pfg_\lambda}
\newcommand{\fn}{\mathbf{0}}
\newcommand{\fv}{\mathbf{v}}
\newcommand{\fd}{\mathbf{d}}
\newcommand{\fr}{\mathbf{r}}

\newcommand{\ful}{\fu^l}
\newcommand{\fuc}{\fu^c}
\newcommand{\fusl}{\fu^{*,l}}
\newcommand{\fusc}{\fu^{*,c}}
\newcommand{\fuslc}{\fu^{*,l/c}}
\newcommand{\fusp}{\fu^{*,+}}
\newcommand{\fusm}{\fu^{*,-}}
\newcommand{\fuspm}{\fu^{*,\pm}}

\newcommand{\fuspT}{\fu^{*,+}_T}
\newcommand{\fusmn}{\fu^{*,-}_0}

\newcommand{\fusln}{\fu^{*,l}_0}
\newcommand{\fuscn}{\fu^{*,c}_0}
\newcommand{\fuslcn}{\fu^{*,l/c}_{0}}

\newcommand{\fuslT}{\fu^{*,l}_T}
\newcommand{\fuscT}{\fu^{*,c}_T}
\newcommand{\fulT}{\ful_T}
\newcommand{\fucT}{\fuc_T}
\newcommand{\fud}{\hat{\fu}}

\newcommand{\fuT}{\fu_T}
\newcommand{\fvn}{\fv_0}

\newcommand{\cC}{\mathcal{C}}
\newcommand{\cD}{\mathcal{D}}
\newcommand{\omr}{\om_{\varrho}}
\newcommand{\Sr}{S_{\varrho}}
\def\issn{1456-436X}